\documentclass[prb,preprintnumbers,amsmath,amssymb,twosides]{revtex4}

\usepackage[applemac]{inputenc}
\usepackage{amsfonts}
\usepackage{graphicx}
\usepackage{dcolumn}
\usepackage{bm}
\usepackage[english]{babel}
\usepackage{amsmath}

\newcommand{\trm}{\textrm}

\newcommand{\be}{\beqn}
\newcommand{\ee}{\eeqn}
\newcommand{\bea}{\begin{eqnarray}}
\newcommand{\eea}{\end{eqnarray}}

\newcommand{\beqn}{\begin{eqnarray}}
\newcommand{\eeqn}{\end{eqnarray}}
\newcommand{\nn}{\nonumber \\}
\newcommand{\ir}{\mathbb{IR}}

\newcommand{\psia}{$\psi$-algorithm }

\newcommand{\mRz}{$\mathcal{R}_0$}
\newcommand{\mPz}{$\mathcal{P}_0$}
\newcommand{\mR}{$\mathcal{R}$}
\newcommand{\mP}{$\mathcal{P}$}

\usepackage{graphicx}

\begin{document}

\title{System adjustment for targeted performance combining symbolic regression and set inversion}

\author{Abdel KENOUFI}
\email{kenoufi@s-core.fr, abdel.kenoufi@uha.fr}
\affiliation{SCORE, Scientific Consulting for Research and Engineering, LAMIE, Universit\'e de Haute-Alsace, Mulhouse, France}

\author{Jean-Fran\c cois Osselin}
\email{jf.osselin@uha.fr}
\affiliation{LPMT, Laboratoire de Physique et M\'ecanique Textiles, Universit\'e de Haute-Alsace, Mulhouse, France}

\author{Bernard Durand}
\email{bernard.durand@uha.fr}
\affiliation{LPMT, Laboratoire de Physique et M\'ecanique Textiles, Universit\'e de Haute-Alsace, Mulhouse, France}

\date{\today}

\begin{abstract}
One presents methodology and algorithms to prepare a causal system in order to achieve desired 
performances if only input-output data are known and when no other informations are available.
This can be done with mean of evolutionnary programming and set inversion methods, such as $\psi$-algorithm or \emph{\emph{SIvIA}}.
\end{abstract}

\maketitle

\section{Introduction}
The Holy Grail of Science and Engineering is the ability to perform adjustments of a system in order to get the desired performances. For example, how to settle a robot in order to make it moving toward the desired target, or which materials properties are well-suited for a defined functionnality. Thus, many problems are  dealing with the inversion of the relation between adjustments and performances.\\
 Some powerful methods based on interval arithmetic have been developed those last years, such as \emph{\emph{\emph{SIvIA}}}\cite{Jaulin} (Set Inversion \emph{Via} Interval Arithmetic) and $\psi$-algorithm\cite{psi} (Probabilistic Set Inversion) which is inspired from \emph{\emph{\emph{SIvIA}}}. They are both based on interval arithmetic\cite{Warmus1,Warmus2,Sunaga,Markov1,Markov2,Markov3,Markov4,Moore1,Moore2,Moore3,Moore4,Moore5,nicolas,rc1}, which has been introduced originally to compute very quickly and in a reliable way quantities represented with intervals. Nowadays, interval arithmetic has been extended to become a computing framework which permits to perform error analysis by computing mathematic bounds, and to solve problems such as non linear problems, partial differential equations, and set inversion\cite{Jaulin,psi}. It finds a large place of applications in material sciences, controllability, automatism, robotics, embedded systems, biomedical, haptic interfaces, form optimization, analysis of architecture plans, ...\\
We are interested to compute the adjustments needed by a system to achieve specific performances, when the modelization of the system is unknown and when only experimental and measured data are available. In order to perform such kind of difficult set inversion, one proposes in this paper to combine \emph{\emph{\emph{SIvIA}}}-like algorithms such as $\psi$-algorithm with some evolutionnary schemes, in particular symbolic regression using genetic programming. In the next section, one recalls briefly the basis of the $\psi$-algorithm, and exhibit its limitations. Thereafter, one presents symbolic regression based on genetic programming to simulate the system that we would like to invert. We exhibit in the third section some examples with numerical applications developed in \emph{ruby}\cite{ruby} which is an elegant, powerful object oriented programming language.

\section{Methodology}

\begin{figure}[!h]
\centering
\includegraphics[scale=1.5]{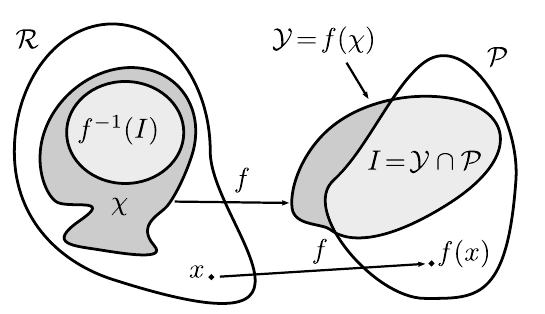}
\caption{Set Inversion consists to compute  $\mathcal{S}=f^{-1}(\mathcal{P})\cap\mathcal{R}$ for a defined function $f$ with mean of interval arithmetics.}
\label{f1}
\end{figure}

Let us note $\mathcal{R}\subset\mathbb{R}^n$ the set of adjustments, and $\mathcal{P}\subset\mathbb{R}^p$ the set of performance of a system. The mathematical modelization of our problem consists of the computation of $\mathcal{S}=f^{-1}(\mathcal{P})\cap\mathcal{R}$, as shown on figure (\ref{f1}), where $f : \mathbb{R}^n\rightarrow \mathbb{R}^p$ is the function giving performances from adjustments. For semantical reasons, one has to perform this set inversion within the interval semi-group $\ir$\cite{psi}. But unfortunately, in many cases, this function is not known, i.e no formulas can be given to reproduce the output values from the input ones. However some input-output measurements, meaning adjustments-performances experimental or empirical values can be available and that's the reason one proposes to use first the symbolic regression based on genetic programming to approximate this function (figure \ref{f2}) in the sense of uniform convergence on the adjustments set, and it will be used afterward within a set inversion scheme such as pure \emph{\emph{\emph{SIvIA}}} or the $\psi$-algorithm (figure \ref{f3}).

\section{Symbolic Regression \emph{via} Genetic Programming (SR\emph{v}GP)}

\subsection{Scheme}

\begin{figure}[!h]
\centering
\includegraphics[scale=1.5]{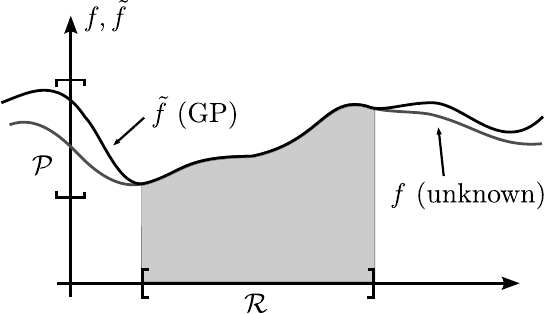}
\caption{Symbolic Regression \emph{via} Genetic Programming permits to build an approximation of an unknown function defined on \mR, and valued on \mP.}
\label{f2}
\end{figure}

Usually if some input-output data are known, one performs function parameters fitting using optimization methods. Because it needs some knowledge of an empirical and parametrized form of the function, and since it is hard for the numerical scheme not to be trapped in a local minimum, it can be rather useful and powerful to perform a symbolic regression in order to find a function converging uniformly to the unknown one as shown on figure (\ref{f2}). \\
Symbolic regression is a good candidate to achieve such kind of procedures. It has the great advantage to give a formal expression of the wanted function which converges uniformely on the adjustment set (figure \ref{f2}) to the unknown one. Uniform convergence means that it is stronger than pointwise convergence and that the speed of convergence is not point dependant. \\
Symbolic regression can be done in different ways. In this article, we focus on the so-called Symbolic Regression \emph{via} Genetic Programming (GP)\cite{Koza2005,Langdon2002,Brownlee}(\emph{SRvGP}). GP was initially  developed for inductive automatic programming and is well-suited for symbolic regression, controller design, and machine learning tasks. One can consider GP algorithms as an extension of the well-known genetics algorithms. GP's scope is to use induction to devise a computer program. This is achieved by using evolutionary operators on candidate programs with a tree structure to improve the adaptive fit between the population of candidate programs and an objective function. An assessment of a candidate solution involves its execution. Symbolic regression with input-output data is achieved just by considering a function as a computer program, which is represented with a tree (figure \ref{arbre}) and to use those data in the objective function.

\begin{figure}[!h]
\centering
\includegraphics[scale=2.0]{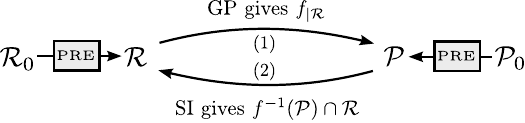}
\caption{Combining Genetic Programming (GP) and Set Inversion (SI) permits first to find an approximation of the model and in a second step to invert performance set in order to find the efficient and possible adjustments. Experimental input-output data could be filtered and decimated if it is needed in order to remove some frequencies and to avoid overfitting phenomena.}
\label{f3}
\end{figure}

Let us suppose that some input-output data (\mRz,\mPz) have been measured. Therefore, we would like to find out a symbolic formula which can describe the system by mean of the symbolic regression. However, it can be useful sometimes to perform pre-treatements with Fourier or wavelets transforms \cite{Altaisky}(PRE-treatment on figure (\ref{f3}). This signal processing needs that one has to choose the frequencies which have to been kept or removed. On the other hand, it can have meaning to find a solution which takes into account the noise, but this is another approach which is not discussed in this paper and can be studied in another framework. In a second step, a decimation of the signal points can be useful to decrease the computational time and to avoid overfitting phenomena. Wavelet analysis seems to be a promising tool\cite{Altaisky} to achieve that. Finally, one obtains a new input-output data set which can be used for the symbolic regression of the system (figure \ref{f3}).

\subsection{Numerical example of SR\emph{v}GP}

For example, we propose to illustrate this procedure on data produced with the function 
\begin{equation}
f(x)=\sin(5x)\cdot e^{-x^2} \ \trm{for } x\in [-3,3]
\end{equation}
This signal is very interesting, not only because it is a wavelet\cite{Altaisky}, but because it has several local extrema and fast oscillations.\\
We have added to this function a random uniformly distributed noise $\epsilon(x) \in [-0.25,0.25]\ \trm{for } x \in [-3,3]$ in order to see if, in this case, the genetic programming scheme is sensible to small perturbations or could extract the right signal we are looking for.

\begin{figure}[!h]
\centering
\includegraphics[scale=0.8]{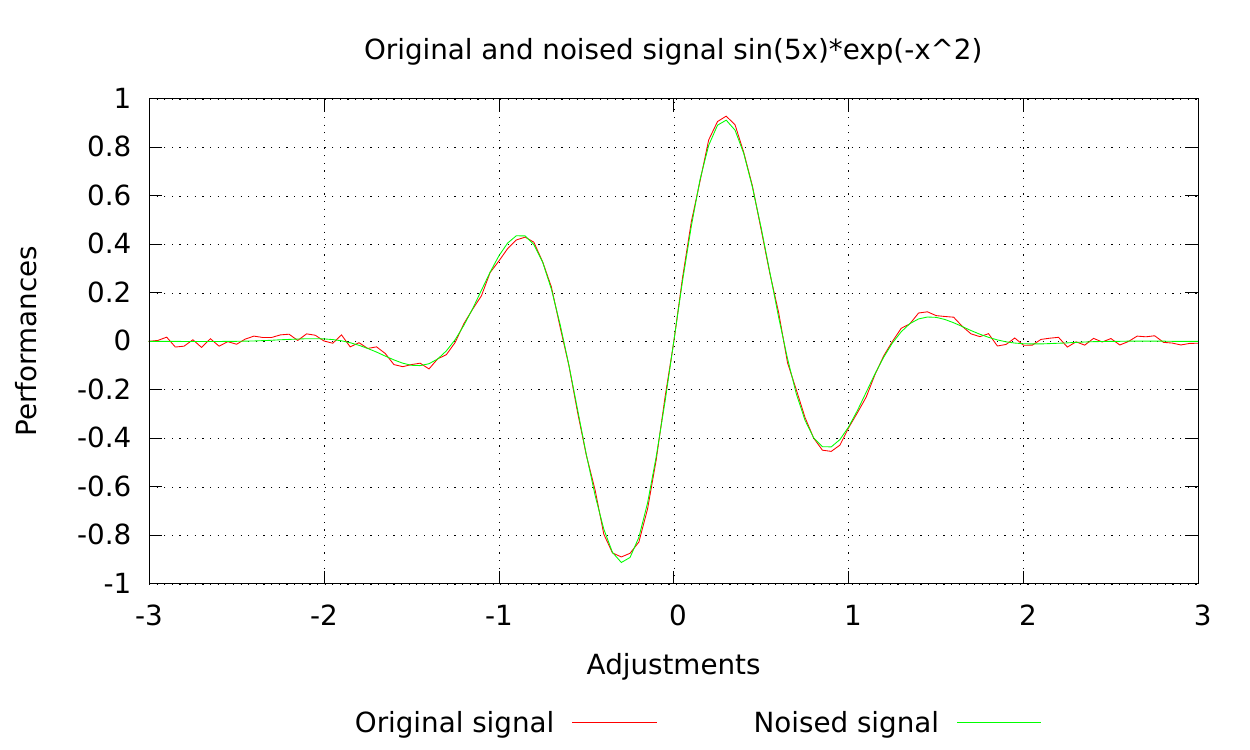}
\caption{Performance set for the function $f(x)=\sin(5\cdot x)\cdot e^{-x^2}$. An uniform random noise $\epsilon(x) \in [-0.25,0.25]\ \trm{for } x \in [-3,3]$ has been added to the initial output data set f([-3,3]).}.
\label{perf}
\end{figure}

\begin{figure}[!h]
\centering
\includegraphics[scale=.5]{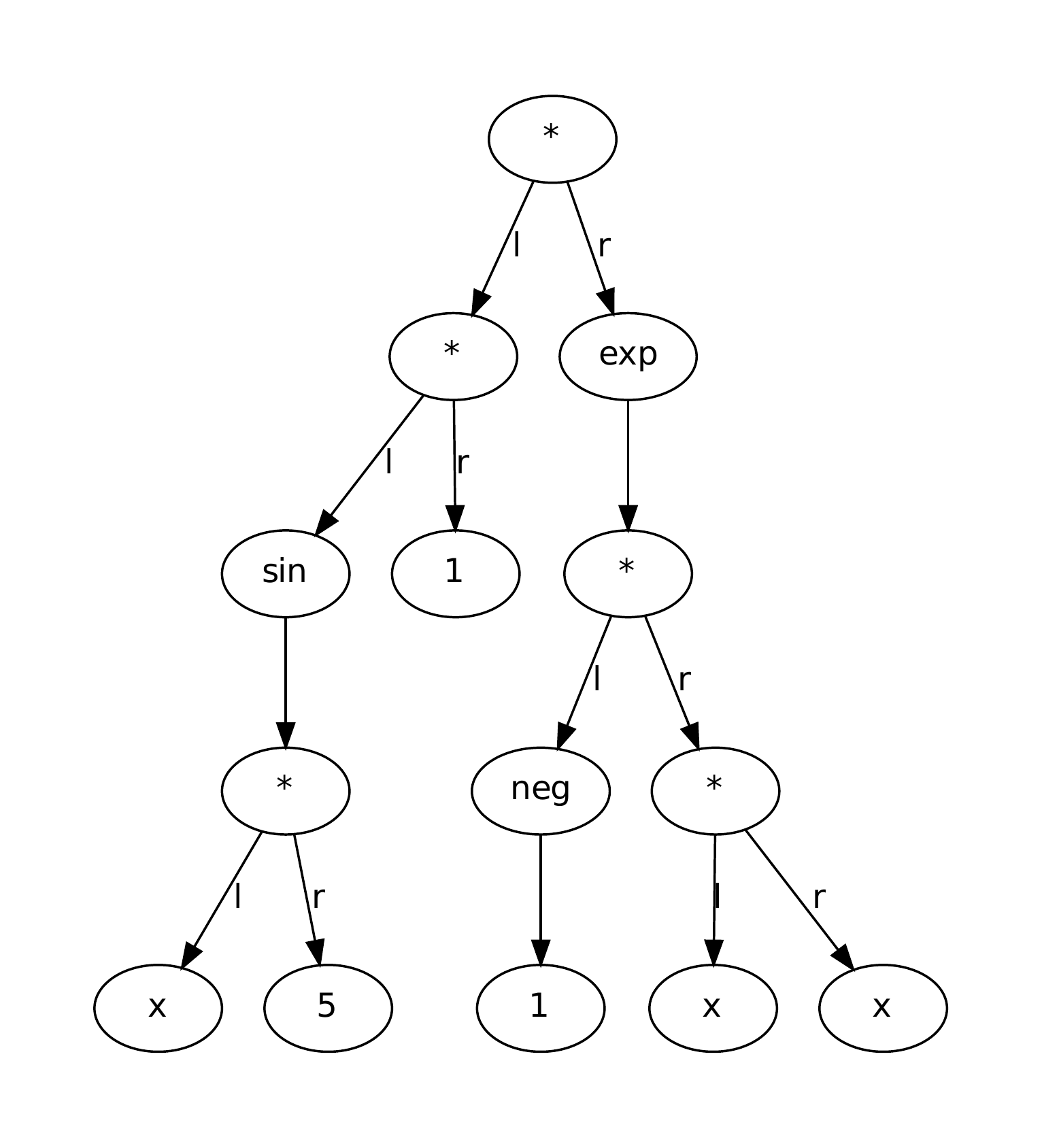}
\caption{Result for the genetic programming of the $f(x)=\sin(5 x)\cdot e^{-x^2}$ noised example. This tree corresponds to the approximation of $f$, $\tilde{f}(x)=\sin(5 x)\cdot e^{-x^2} = f(x)$ on $[-3,3]$ with function basis set $\{+, -, \cdot, /, \exp, \log, \sin, \cos, \tan, \trm{neg}\}$. The letters "l" and "r" correspond respectively to left and right positions of operands in a binary expression.}
\label{arbre}
\end{figure}

The symbolic regression has been implemented within a \emph{ruby} script and a MPI library for this language that we have developed for the scope of this present article. The code is well-suited for parallel computer, launching several jobs on different nodes with different random seeds. After convergence of each population, a last run is done on a single node involving the best candidates of each previous nodes within a larger population. \\
The final result obtained in less than 3 minutes for this example with 16 processors and for a cost function about $10^{-3}$, is given in a \emph{EBNF}\cite{EBNF} (Extended Backus-Naur Form) expression $(*\  (*\  (\sin \ (*\  x\  5))\  1)(\exp \  (* \  (\trm{neg}\  1)(*\  x\  x))))$ and can be represented by the tree shown on figure (\ref{arbre}). This expression represents the function 
\begin{equation} \tilde{f}(x)=\sin(5x)\cdot e^{-x^2}.
\label{wavelet}
\end{equation}
We have choosen a grow method for the tree developpement with a maximal depth of 5 and a set of function such as $\{+, -, \cdot, /, \exp, \log, \sin, \cos, \tan, \trm{neg}\}$ where $\trm{neg}\ : x\mapsto -x$. Figure (\ref{perf}) shows the noised signals for adjustments on $[-3,3]$. It is interesting that in this example the signal filtering has not been used in this example since the genetic program absorbs the noise and finds the right solution.  Thus, the symbolic regression finds finally the right function $\tilde{f} = f$ on the adjustments set $[-3,3]$ even in presence of small noise and without filtering. 

\section{Probabilistic Set Inversion}
\subsection{Scheme}
Let us assume now that the model describing the adjustment-performance is known. In order to perform set inversion we propose to use the \psia\cite{psi} (Probabilistic Set Inversion) inspired from \emph{\emph{SIvIA}}\cite{Jaulin} but without boolean tests and with a conditional probability calculation and domain bissections. The interval arithmetics used is based on free algebras\cite{nicolas,rc1} and permits to define univoquely inclusion function associated to the natural one and to avoid data dependencies\cite{nicolas,rc1,psi}. This yields to accepted or rejected domains only. In this algorithm, one computes the following conditional probability with notations related to the figure (\ref{f1}), 

\beqn
p(\mathcal{X})&=&p([f](\mathcal{X})\subset\mathcal{P}\ | \  f(x)\in[f](\mathcal{X}),\ \forall x\in \mathcal{X})\nn
&=&\frac{mes([f](\mathcal{X})\cap\mathcal{P})}{mes([f](\mathcal{X}))}\nn
&=&\frac{mes(\mathcal{Y}\cap\mathcal{P})}{mes(\mathcal{Y})}=\frac{mes(\mathcal{I})}{mes(\mathcal{Y})}.\nn
\eeqn

where $mes$ is the Lebesgue measure in $\mathbb{R}^p$ (p-dimension volume).\\
 If the image of an adjustement subset $\mathcal{X}$ is included in the performance set $\mathcal{P}$, then $\mathcal{I}=\mathcal{Y}$. This yields the probability to be equal to $1$ and one adds the subset $\mathcal{X}$ to the solution list. If the image set $\mathcal{Y}$ is not included in the performance set, then $mes(\mathcal{I})=0$ and the probability equals zero, therefore the subset $\mathcal{X}$ is rejected and removed from the list of interval candidates. If $p(\mathcal{X}) \in\  ]0,1[$, then $\mathcal{X}$ is bissected in all directions and \psia applies the same procedure recursively for the resulting intervals until the size is lower than a fixed volume resolution of the intervals or until the sets are accepted or rejected. Moreover, \psia uses a well-defined interval arithmetic based on free algebra\cite{nicolas,rc1,psi}. It guarantees to get well-defined inclusion functions and avoids in particular data dependancies\cite{nicolas,rc1,psi,Jaulin}.\\\

\subsection{Numerical examples of \psia  applications}

\subsubsection{One dimensional example}

Let's suppose now that the model of our problem is the previous function found in the equation (\ref{wavelet}) with a SR\emph{v}GP scheme. Now one would like to solve the following inequation :

\begin{equation}
-\frac{1}{4}  \le  \tilde f(x)  \le  \frac{1}{2},\ \trm{for } x\in [-2,2]
\end{equation}

This can be rewritten as a set inversion for $\tilde f$ with $\mathcal{R}=[-2,2]$ and $\mathcal{P}=[-\frac{1}{4},\frac{1}{2}]$.

\begin{figure}[!h]
\centering
\includegraphics[scale=.8]{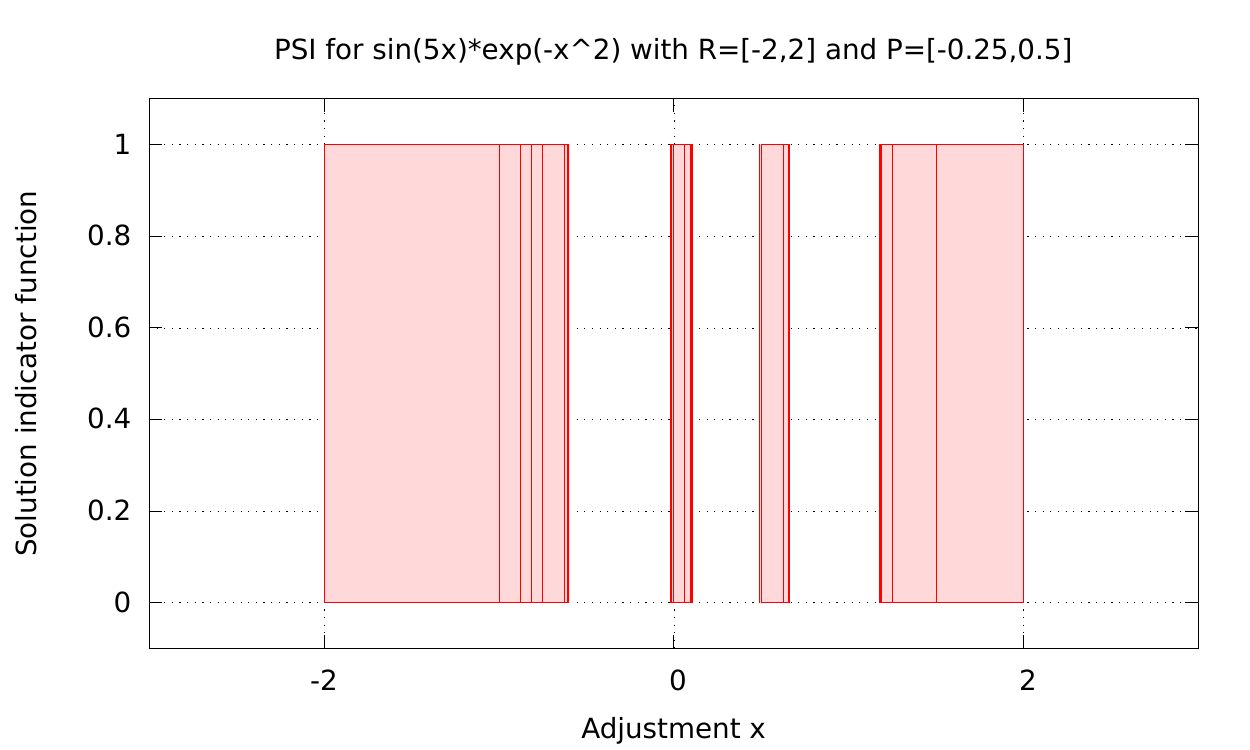}
\caption{Result of set inversion with \psia for $\tilde f (x)=\sin(5 x)\cdot e^{-x^2}$ with $\mathcal{R}=[-2,2]$ and $\mathcal{P}=[-\frac{1}{4},\frac{1}{2}]$}
\label{psi_wavelet}
\end{figure}

\subsubsection{Higher dimensional examples}

Let us give an other example with the function $f : \mathbb{R}^2\longrightarrow \mathbb{R}^2$ : 
\be
f(x,y)=(x^2-y^2\cdot e^x+x\cdot e^y,x\cdot(x+y)-y^2)
\ee 

\begin{figure}[!h]
\centering
\includegraphics[scale=.8]{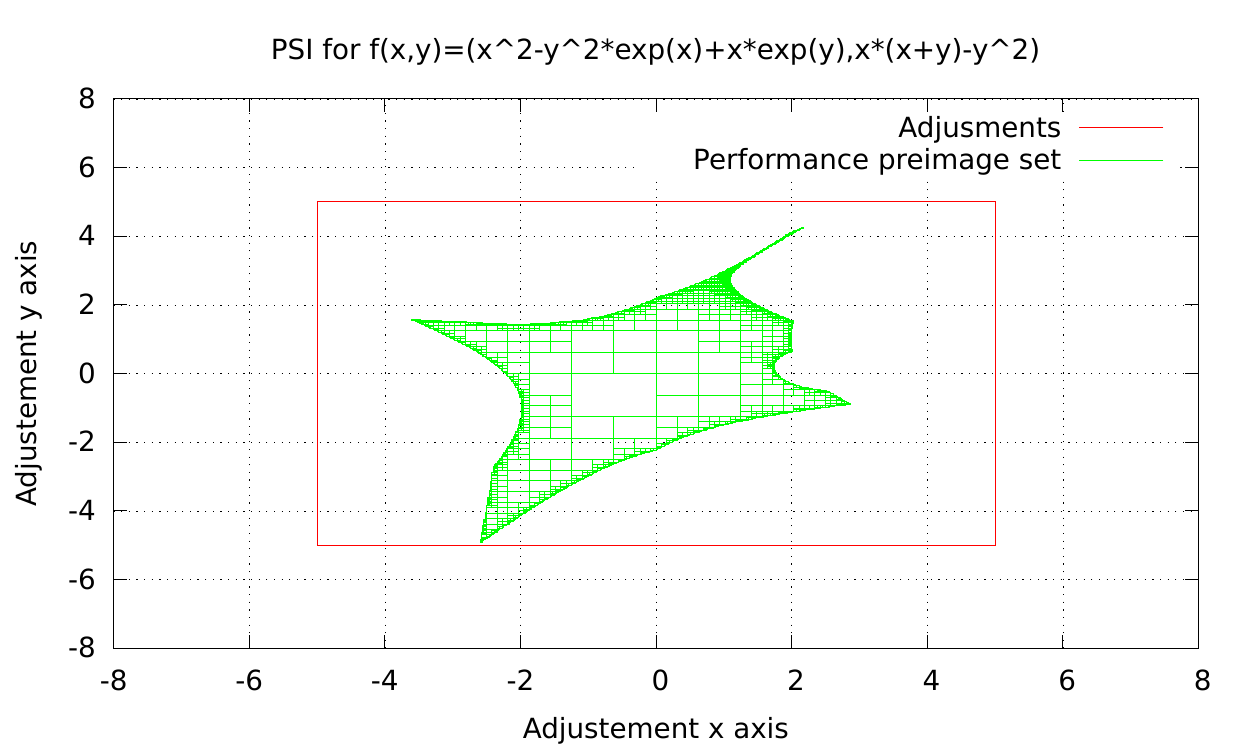}
\caption{\psia for $f(x,y)=(x^2-y^2\cdot e^x+x\cdot e^y,x\cdot(x+y)-y^2)$ with $\mathcal{R}=\mathcal{P}=[-5,5]^2$.}
\label{psi4}
\end{figure}

One would like to invert this function for $\mathcal{R}=[-5,5]^2$ and $\mathcal{P}=[-5,5]^2$. This corresponds to the resolution of a strongly non-linear inequations system. The set inversion presented on figure (\ref{psi4}) shows clearly that the \psia implemented in the algebraic arithmetic we use is not data dependant. The variables appear more than once in the formal expression of the function $f$. This shows how \emph{SIvIA}-inspired methods such as \psia are powerful and efficient to inverse performance set through a continuous function. The CPU time for this inversion is about $130$ seconds on a two processors computer ($1.67$ Ghz Intel) for a spatial surface resolution of $10^{-4}$.

\begin{figure}[!h]
\centering
\includegraphics[scale=.8]{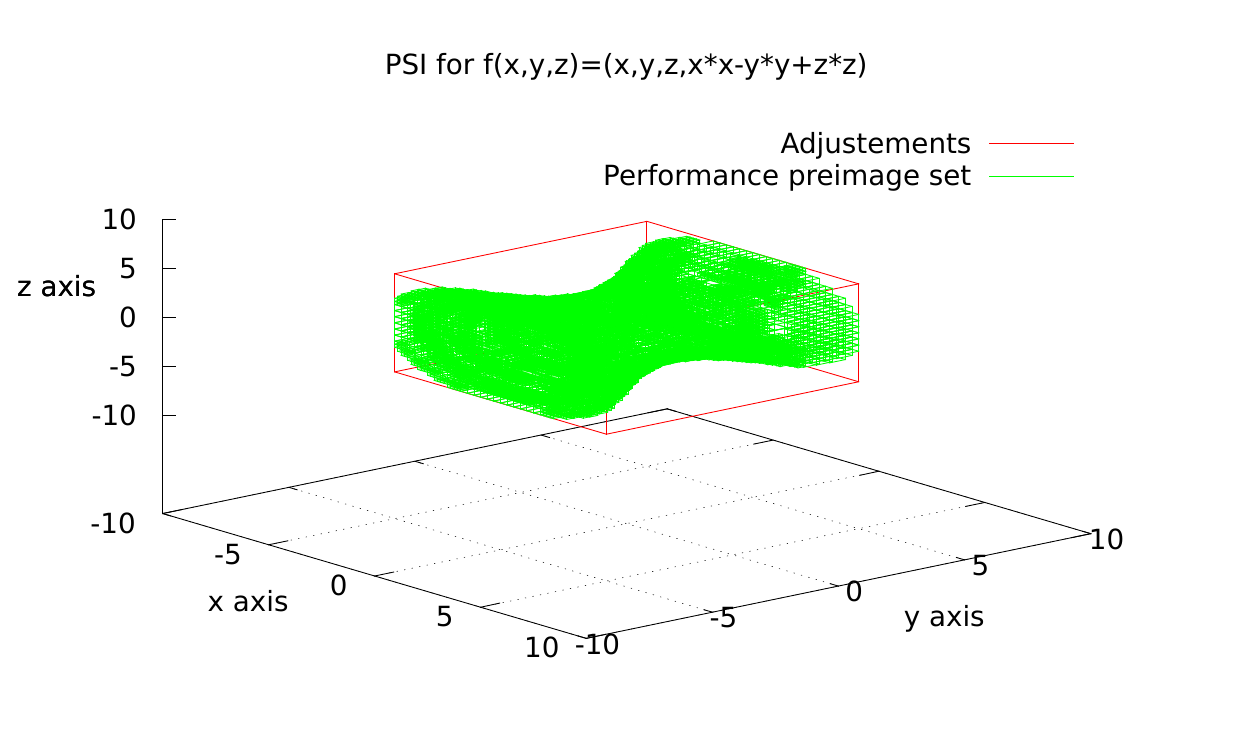}
\caption{\psia for $f(x,y,z)=(x,y,z,x^2-y^2+z^2)$ with $\mathcal{R}=[-5,5]^3$ and $\mathcal{P}=[-10,10]^4$.}
\label{psi6}
\end{figure}

There is no limitation for the dimensions of the adjustments and performances sets as shown on figure (\ref{psi6}) for the function $f : \mathbb{R}^3\longrightarrow \mathbb{R}^4$. Due to the bissection of the domain, the computational time remains still exponential. However, it is possible to decrease the computational time constant with mean of parallelization using domain decomposition\cite{Jaulin,psi}. The adjustement set is divided on the first axis, and each processor performs the \psia on one of those subdomains. The master processor collects all the results at the end of the calculations.

\section{Conclusion and perspectives}

    The paper has presented a new methodology for adjustments of a system in order to satisfy desired performances. This technique, combining genetic programming for symbolic regression with set inversion, only assumes that input-output data are known. The approach is decomposed into two steps: the symbolic regression and the set inversion. The role of the symbolic regression is to find an accurate model for the system by fitting all input-output data. The main advantage of symbolic regression compared to more classical fitting techniques such as splines, wavelet, etc. is that the resulting model has a short symbolic expression which is essential to make interval set inversion method efficient. The set inversion with \psia is able to solve many interesting problems occuring in science and engineering such as strongly non-linear inequations and equations, and optimization problems.
    
\begin{acknowledgements}
We thank Luc Jaulin from ENSTA-Bretage, Michel Goze from Universit\'e de Haute-Alsace, Michel Gondran and Thierry Socroun from Electricit\'e De France for useful and interesting discussions.
\end{acknowledgements}


\end{document}